\documentclass[12pt,a4paper,makeidx,reqno]{amsart}

\makeindex \topmargin=0in \headheight=8pt \headsep=0.5in \textheight=8.9in
\usepackage{amsmath,amssymb,amsthm}
\usepackage{mathrsfs}
\usepackage{eucal}
\usepackage[matrix,arrow,curve]{xy}
\usepackage{textcomp}
\usepackage{graphicx}
\usepackage{wrapfig}
\usepackage{verbatim}
\usepackage{tikz-cd}

\theoremstyle{remark}

\numberwithin{equation}{section}

\def\Z{{\mathbb Z}} 
 \def\Q{{\mathbb Q}} \def\A{{\mathbb A}}

\long\def\comment#1\endcomment{}

\def\Z{{\mathbb Z}} 
 \def\Q{{\mathbb Q}} \def\A{{\mathbb A}}

\begin{document}

\title{ \bf Homologies of inverse limits of groups}

\author{Danil Akhtiamov }

\footnotetext{The author is supported by the grant of the Government of the Russian Federation
for the state support of scientific research carried out under the supervision of leading scientists,
agreement 14.W03.31.0030 dated 15.02.2018.}

\smallskip

\maketitle

\begin{abstract}
   Let $H_n$ be the $n$-th group homology functor (with integer coeffcients) and let $\{G_i\} _ {i \in \mathbb{N}}$ be any tower of groups such that all maps $G_{i+1} \to G_i$ are surjective.
   In this work we study kernel and cokernel of the following natural map: $$H_n(\varprojlim G_i) \to \varprojlim H_n(G_i)$$ For $n=1$ Barnea and Shelah [BS] proved that this map is surjective and its kernel is a cotorsion group for any such tower $\{G_i\} _ {i \in \mathbb{N}}$. We show that for $n=2$ the kernel can be non-cotorsion group even in the case when all $G_i$ are abelian   and after it we study these kernels and cokernels for towers of abelian groups in more detail.
\end{abstract}
 
\section{Introduction}

It is well-known that $H_n$ commute with direct limits for any $n \ge 0$, where $H_n=H_n(-, \Z)$ s a functor from topological space or from groups to abelian groups.  But in general it is not true for projective limits. Moreover, it is quite difficult to understand whether they commute or not in some concrete cases. Some of these cases are connected with difficult problems. For example, it is known that $\varprojlim{H_3(F/\gamma_n(F))}=0$ for finitely generated free group $F$. Although it turned out difficult to understand whether $H_3(\varprojlim{F/\gamma_n(F)})=0$ or not and it is still an open problem. And  $H_3(\varprojlim{F/\gamma_n(F)}) \ne 0$ will imply that answer to the Strong Parafree Conjecture, which is also still open, is negative [Hill, p. 294]. Also in some cases it is quite easy to see that kernel of the natural map $H_n(\varprojlim{G_i}) \to \varprojlim{H_n(G_i)}$ is non-zero, but it is hard to say something about structure of this kernel. For example, let $n=1$ and $G_i=F^{\times i}$, where, again, $F$ is a finitely generated free group. It is proven by Miasnikov and Kharlampovich [KM] that the kernel in this case contains $2$-torsion but their prove is quite difficult and uses very heavy machinery called "Non-commutative Implicit Function Theorem". And, additionally, question about $p$-torsions for prime $p>2$ is still open. 

This paper is an attempt to start systematic study of kernels and cokernels of maps $\mathfrak{F}(\varprojlim{G_i}) \to \varprojlim{\mathfrak{F}(G_i)}$, where $F$ is a functor, especially in the case $\mathfrak{F}=H_n$. It is well-known that for functor $\pi_n: Top_{*} \to Ab$ and for a tower of connected pointed spaces $X_i$, such that maps $X_{i+1} \to X_i$ are fibrations there is the following short exact sequence, called Milnor sequence (see [GJ, VI Proposition 2.15]): $$0 \to {\varprojlim}^1{\pi_{n+1}(X_i)} \to \pi_{n}(\varprojlim{X_i}) \to  \varprojlim{\pi_{n}(X_i)} \to 0$$
The similar fact is true for the functor $H_n: Ch_{\mathbb{Z}} \to Ab$, where $Ch_{\mathbb{Z}}$ is the category of chain complexes of abelian groups and $H_n$ is $n-th$ chain homology functor (see, for instance, [W, p. 94, prop. 3.5.8]): Let $C_i$ be a tower of chain complexes (of abelian groups) such that it satisfies degree-wise the Mittag-Leffler condition. Then there is the following short exact sequence:
$$0 \to {\varprojlim}^1{H_{n+1}(C_i)} \to H_{n}(\varprojlim{C_i}) \to  \varprojlim{H_{n}(C_i)} \to 0$$

So, it is natural to ask if the same statement true for  homologies of spaces with integer coefficients  for a tower of connected pointed spacs $X_i$, such that maps $X_{i+1} \to X_i$ are fibrations. (Un)fortunately, the answer is no for two reasons. First, the kernel might be non-zero, while  ${\varprojlim}^1{H_{n+1}(X_i)}=0$. To see it let us consider $X_i=K(F/\gamma_i(F))$, where $F$ is $2$-generated free group. Using fibrant replacments, we can assume that all maps $X_{i+1} \to X_{i}$ are fibratons. Then, using Milnor sequences, we see that $\varprojlim{X_{i}}=K(\varprojlim{F/\gamma_n(F)},1)$. This case was studied by Bousfield and Kan deeply because of its connection with Bousfield's completion of spaces. It is known that ${\varprojlim}^1{H_{n+1}(X_i)}=0$ for every $n \ge 0$ [BK, p.123]. But it is proven by Bousfied in [Bous] that $H_2(F/\gamma_n(F)) \ne 0$ (and it implies that wedge of two circles is $\mathbb{Z}$-bad space; actually wedge of two circles is also  $\mathbb{Q}$-bad and $\mathbb{Z}/p$-bad space for $p>2$, but it was proven much more later by Ivanov and Mkhailov in the papers [IM2], [IM3]). Second, the cokernel might be non-zero. Corresponding example actually was provided by Dwyer in [D] (Example 3.6). See Corollary 4 from my work for more details on this example.

{\bf Definition. } We call an abelian group $A$ a cotorsion group if $A= {\varprojlim}^1{B_i}$ for some tower of abelian groups $B_i$.

{\it Remark. } Usually people define cotorsion groups in a different way but these definitions coincide because of [H] (Theorem 1). More detailed, there is the following equivalent definition of cotorsion groups: 

{\bf Definition. } We call an abelian group $A$ a cotorsion group if $Ext(C,A)$=0 for any torsion-free abelian group $C$.

There is the following result which was proven by Shelah and Barnea  ([BS], Corollary 0.0.9):

{\bf Theorem. } {\it  Let $X_i$ be an inverse system  of pointed connected spaces, such that all maps $X_{i+1} \to X_i$ are Serre fibrations, $\pi_1(X_{i})$ and $\pi_2(X_{i})$ satisfies the Mittag-Leffler condition. Then the natural map $H_1(\varprojlim{X_i}) \to \varprojlim {H_1(X_i)}$ is surjective and its kernel is a cotorsion group.}

This result might give one a hope that, assuming towers $\pi_1(X_{i}), \pi_2(X_{i}), \dots,  \pi_k(X_{i})$ for large enough (maybe infinite) $k$ satisfy Mittag-Leffler condition, we can "fix" usual Milnor sequences and provide "Milnor sequences for $H_n$". 

(Un)fortunately, the answer is no already for $n=2$. Moreover, it is false even for homologies of abelian groups and it is proven in this work:

{\bf Theorem 1. }{\it There is an inverse system of abelian groups indexed by $\mathbb{N}$ such that all maps $A_{i+1} \to A_i$ are epimorhisms and kernel of the natural map $H_2(\varprojlim{A_i}) \to \varprojlim{H_2(A_i)}$ is not a cotorsion group. }

This theorem shows that things are very difficult even for abelian groups, so the following result seems quite interesting:

{\bf Theorem 2. }{\it Let $A_i$ be an inverse system of torsion-free abelian groups indexed by $\mathbb{N}$ such that all maps $A_{i+1} \to A_i$ are epimorhisms. Then for any $n \in \mathbb{N}$ the natural map $H_n(\varprojlim{A_i}) \to \varprojlim{H_n(A_i)}$ is an embedding and its cokernel is a cotorsion group. }

\smallskip

{\bf Theorem 3. }{\it Let $A_i$ be an inverse system of any abelian groups indexed by $\mathbb{N}$ such that all maps $A_{i+1} \to A_i$ are epimorhisms. Then:

(1) Cokernel of the natural map $H_2(\varprojlim{A_i}) \to \varprojlim{H_2(A_i)}$ is a cotorsion group

(2) Suppose, additionally, that torsion subgroup of $A_i$ is a group of bounded exponent for any $i$. Then cokernel of the natural map $H_3(\varprojlim{A_i}) \to \varprojlim{H_3(A_i)}$ is a cotorsion group. In particular, it is true if all $A_i$ are finitely generated.  }

Also in this work there are two results not about homologies but about another functors and inverse limits of abelian groups:

{\bf Statement 2. }{\it Let $B$ be any abelian group and $B_i$ be an inverse system of abelian groups with surjective maps $f_i: B_{i+1} \to B_i$ between them. Then kernel of the map $Tor(B,\varprojlim{B_i}) \to \varprojlim{Tor(B,B_i)}$ is trivial. }

{\bf Corollary 6. }{\it $Tor (\varprojlim{A_i},\varprojlim{A_i}) \to \varprojlim {Tor(A_i,A_i)}$ is embedding for any inverse system of abelian groups such that all maps $A_{i+1} \to A_i$ are surjective.}

\smallskip

The author thanks Roman Mikhailov, Sergei O. Ivanov, Emmanuel Farjoun, Fedor Pavutnitsky and Saharon Shelah for fruitful discussions.

\section{Homologies of inverse limits of groups}

 A goal of this chapter is to prove the following results:

{\bf Theorem 1. }{\it There is an inverse system of abelian groups indexed by $\mathbb{N}$ such that all maps $A_{i+1} \to A_i$ are epimorhisms and kernel of the natural map $H_2(\varprojlim{A_i}) \to \varprojlim{H_2(A_i)}$ is not a cotorsion group. }

\smallskip

{\bf Theorem 2. }{\it Let $A_i$ be an inverse system of torsion-free abelian groups indexed by $\mathbb{N}$ such that all maps $A_{i+1} \to A_i$ are epimorhisms. Then for any $n \in \mathbb{N}$ the natural map $H_n(\varprojlim{A_i}) \to \varprojlim{H_n(A_i)}$ is embedding and its cokernel is a cotorsion group. }

\smallskip

{\bf Theorem 3. }{\it Let $A_i$ be an inverse system of any abelian groups indexed by $\mathbb{N}$ such that all maps $A_{i+1} \to A_i$ are epimorhisms. Then:

(1) Cokernel of the natural map $H_2(\varprojlim{A_i}) \to \varprojlim{H_2(A_i)}$ is a cotorsion group

(2) Suppose, additionally, that torsion subgroup of $A_i$ is a group of bounded exponent for any $i$. Then cokernel of the natural map $H_3(\varprojlim{A_i}) \to \varprojlim{H_3(A_i)}$ is a cotorsion group. In particular, it is true if all $A_i$ are finitely generated.  }

\smallskip

{\it Proof of Theorem 1. } Let us consider $$\displaystyle A'_{i,p}=\mathbb{\Z}^i \oplus (\bigoplus_{i+1}^{\infty} \mathbb{Z}/ p\mathbb{Z}), A'_i:=\bigoplus_{p \in \mathbb{P}} A'_{i,p}, B:=\bigoplus_{p \in \mathbb{P}} \mathbb{Z}/ p\mathbb{Z}, A_i:=A_i' \times B$$
We are going to define maps $\psi_i: A'_{i} \to A'_{i-1}$ in the most natural way. Let us denote by $e_{i,p}^1, \dots, e_{i,p}^i$ elements of basis of the free abelian summand of $A'_{i,p}$  and let us denote by $e_{i,p}^{i+1}, e_{i,p}^{i+2} , \dots$ elements of basis of the $\mathbb{Z}/ p\mathbb{Z}$-vector space, which is the second direct summand of $A'_{i,p}$. Now let us define $\psi_i(e_{i,p}^{j}):=e_{i-1,p}^{j}$. It is obvious that there is unique $\psi_i$ with such properties. And, finally, let $\phi_i: A_{i} \to A_{i-1}$ be the map defined by the matrix 
$\left[ {\begin{array}{cc}
  \psi_i & 0 \\
   0 & Id_{B} \\
  \end{array} } \right]. $

Using Kunneth formula and usig that all $A'_i$ and $B$ are abelian, we have: $$H_2(\varprojlim{A_i})=H_2(\varprojlim{(A'_i \times B)})=H_2((\varprojlim{A'_i}) \times B)=H_2(\varprojlim{A'_i}) \oplus H_2(B) \oplus ((\varprojlim{A'_i}) \otimes B)$$ 

In the similar way we can get: $$\varprojlim{H_2(A_i)}=\varprojlim{H_2((A'_i \times B))}=\varprojlim{H_2(A'_i) \oplus H_2(B) \oplus (A'_i \otimes B)}=\varprojlim{H_2(A'_i)} \oplus H_2(B) \oplus \varprojlim{(A'_i \otimes B)}$$  

Thus our map $H_2(\varprojlim{A_i}) \to \varprojlim{H_2(A_i)}$ is a map from $H_2(\varprojlim{A'_i}) \oplus H_2(B) \oplus ((\varprojlim{A'_i}) \otimes B)$ to $\varprojlim{H_2(A'_i)} \oplus H_2(B) \oplus \varprojlim{(A'_i \otimes B)}$. Analyzing maps from the Kunneth formula, we see that this map is given by a diagonal matrix, the corresponding maps $H_2(\varprojlim{A'_i}) \to \varprojlim{H_2(A'_i)}$ and $(\varprojlim{A'_i}) \otimes B \to \varprojlim{(A'_i \otimes B)}$ coincide with obvious maps which come from definition of inverse limit and the corresponding map $H_2(B) \to H_2(B)$ is isomorphism. Thus kernel of the natural map 
$(\varprojlim{A'_i}) \otimes B \to \varprojlim{(A'_i \otimes B)}$ is direct summand of kernel of the natural map $H_2(\varprojlim{A_i}) \to \varprojlim{H_2(A_i)}$. So it is enough to prove that kernel of the natural map 
$(\varprojlim{A'_i}) \otimes B \to \varprojlim{(A'_i \otimes B)}$ is not a cotorsion group and now we will do it. Let us note that the map $(\varprojlim{A'_i}) \otimes B \to \varprojlim{(A'_i \otimes B)}$ can be decomposed in the following way:

$$\displaystyle (\varprojlim{A'_i}) \otimes B= \bigoplus_{p \in \mathbb{P}} (\varprojlim{A'_i}) \otimes \mathbb{Z}/p\mathbb{Z}  \to \bigoplus_{p \in \mathbb{P}} \varprojlim({A'_i \otimes \mathbb{Z}/p\mathbb{Z} } ) \to \varprojlim({ \bigoplus_{p \in \mathbb{P}} A'_i \otimes \mathbb{Z}/p\mathbb{Z}})  = \varprojlim{(A'_i \otimes B)}$$
It is easy to see that the map $\displaystyle \bigoplus_{p \in \mathbb{P}} \varprojlim({A'_i \otimes \mathbb{Z}/p\mathbb{Z}}) \to \varprojlim({ \bigoplus_{p \in \mathbb{P}} A'_i \otimes \mathbb{Z}/p\mathbb{Z}})$ is injective. Then kernel of the map $\displaystyle (\varprojlim{A'_i}) \otimes B \to \varprojlim{(A'_i \otimes B)}$ equals $\displaystyle \bigoplus_{p \in \mathbb{P}} Ker [ (\varprojlim{A'_i}) \otimes \mathbb{Z}/p\mathbb{Z} \to \varprojlim{(A'_i \otimes \mathbb{Z}/p\mathbb{Z})}]$. It is proven at [IM] (Corollary 2.5) that   $Ker [(\varprojlim{A'_i}) \otimes \mathbb{Z}/p\mathbb{Z} \to \varprojlim{(A'_i \otimes \mathbb{Z}/p\mathbb{Z})}] \cong {\varprojlim}^1(Tor(A'_i, \mathbb{Z}/p\mathbb{Z}))$. Let us note that $\displaystyle Tor(A'_i, \mathbb{Z}/p\mathbb{Z})=\bigoplus_{i+1}^{\infty} \mathbb{Z}/ p\mathbb{Z}$. So  ${\varprojlim}^1(Tor(A'_i, \mathbb{Z}/p\mathbb{Z})) \cong {\varprojlim}^1 (\bigoplus_{i+1}^{\infty} \mathbb{Z}/ p\mathbb{Z})$, which is, obviously, $p$-torsion abelian group and which is nonzero because of [MP] (p. 330, Proposition A.20). Then, finally, we have that kernel of the map $\displaystyle (\varprojlim{A'_i}) \otimes B \to \varprojlim{(A'_i \otimes B)}$ is not a cotorsion group because of [B] (Theorem 8.5). Q.E.D.
\smallskip

We will need the following lemmas and statements in order to prove Theorem 2 and Theorem 3. All inverse systems supposed to be indexed by $\mathbb{N}$.

{\bf Statement 1. }{\it Let $C$ be any category with projective limits and let $B_i$ be a tower of objects from $C$. Let $F: C \to Ab$ be a functor such that $\displaystyle \varprojlim_{n} {Coker [F(\varprojlim{B_i}) \to F(B_n)]}=0$. Then cokernel of the natural map $F(\varprojlim{B_i}) \to \varprojlim{F(B_i)}$ is a cotorsion group.}

This statement was proved in the third version of Barnea's and Shelah's preprint [BS] in the case when functor $F$ preserves surjection. However, their proof was quite complicated. Our proof of the statement is straightforward. 

{\it Proof. } Let denote by $\Phi_i:=Ker(F(\varprojlim{B_i}) \to F(B_i))$ and by $\Psi_i:=Im(F(\varprojlim{B_i}) \to F(B_i))$. Then we have the following exact sequences: $$0 \to \Phi_i \to F(\varprojlim{B_i}) \to \Psi_i \to 0$$ and $$0 \to \Psi_i \to F(B_i) \to Coker [F(\varprojlim{B_i}) \to F(B_i)] \to 0.$$

From the second sequence we get by the assumption $\displaystyle \varprojlim_{n} {Coker [F(\varprojlim{B_i}) \to F(B_n)]}=0$ $\varprojlim{\Psi_i}=\varprojlim{F(B_i)}$. From the first we deduce a new exact sequence: $$0 \to \varprojlim{\Phi_i} \to F(\varprojlim{B_i}) \to \varprojlim{\Psi_i} \to {\varprojlim}^1{\Phi_i} \to 0$$
Thus finally we have the following exact sequence $$0 \to \varprojlim{\Phi_i} \to F(\varprojlim{B_i}) \to \varprojlim{F(B_i)} \to {\varprojlim}^1{\Phi_i} \to 0$$
Now it is enough to note that $\varprojlim^1$ of any inverse system of abelian groups is a cotorsion group by [H] (Theorem 1). Q.E.D.

{\bf Corollary 1. } {\it Let $B_i$ be an inverse system of groups (respectively abelian groups) and any maps between them. Let $F: Grp \to Ab$ (respectively $F: Ab \to Ab$) be a functor such that $\displaystyle \varprojlim_{n}{Coker [F(\varprojlim{B_i}) \to F(B_n)]}=0$. Then kernel of the map $F(\varprojlim{B_i}) \to \varprojlim{F(B_i)}$ equals $\varprojlim{\Phi_i}$, where $\Phi_n=Ker[F(\varprojlim{B_i}) \to F(B_n)]$.}

{\bf Corollary 2. } {\it Cokernel of the map $\Lambda^n(\varprojlim{B_i}) \to \varprojlim{\Lambda^n(B_i)}$ is a cotorsion group for any inverse system $B_i$, such that $B_{i+1} \to B_i$ are epimorphisms.} 

{\it Proof. } Since $B_{i+1} \to B_i$ are epimorphisms, the maps $\Lambda^n (\varprojlim{B_i}) \to \Lambda^n(B_i)=0$ are also epimorphisms (it easily follows from constructive description of projective limits in the category of groups). Then $\Lambda^n (\varprojlim{B_i}) \to \Lambda^n(B_i)$ is epimorphism, because $\Lambda^n$ is right-exact. Then $Coker[\Lambda^n (\varprojlim{B_i}) \to \Lambda^n(B_i)]=0$ and we are done because of Statement 1. 

{\bf Corollary 3. } {\it Cokernel of the map $H_2(\varprojlim{B_i}) \to \varprojlim{H_2(B_i)}$ is a cotorsion group for any inverse system of abelian groups $B_i$, such that $B_{i+1} \to B_i$ are epimorphisms.}

{\it Proof. } It is well-known [Breen, section 6] that $H_2$ is naturally isomorphic to $\Lambda^2$ in the category of abelian groups, so the corollary follows from Corollary 2.

{\it Definition.} Let $G$ be any group. Let $\gamma_1(G):=G$ and $\gamma_{i+1}(G):=[G,\gamma_{i}(G)]$. Series $\gamma_1(G):=G$ are called lower central series of a group $G$.  Let denote $\widehat{G}:=\varprojlim{G/ \gamma_i(G)}$.  

{\bf Corollary 4. } {\it Cokernel of the map $H_2(\widehat{G}) \to \varprojlim{H_2(G/ \gamma_i(G))}$ is a cotorsion group for any group $G$.} 

{\it Proof}. It is clear that the map $H_2( G) \to H_2(G/ \gamma_i(G))$ factors through the map $H_2( \widehat{G}) \to H_2(G/ \gamma_i(G))$. Then it is enough to prove that $\varprojlim {Coker [H_2( G) \to H_2(G/ \gamma_i(G)]}=0$. Let note that maps between groups $Coker [H_2( G) \to H_2(G/ \gamma_i(G)]$ are zero. Really, using that $H_1(G)=H_1(G/ \gamma_i(G))$, we see from the $5$-term exact sequence [Brown, p.47, exercise 6a] that these cokernels are equal to $\frac{\gamma_i(G) \cap [G,G]}{[\gamma_i(G), G]}=\frac{\gamma_i(G)}{\gamma_{i+1}(G)}$. Q.E.D.

Following statement is not really necessary for a proof of the Theorems, but it is interesting by itself and makes clearer what is happening. 

{\bf Statement 2. }{\it Let $B$ be any abelian group and $B_i$ be an inverse system of abelian groups with surjective maps $f_i: B_{i+1} \to B_i$ between them. Then kernel of the map $Tor(B,\varprojlim{B_i}) \to \varprojlim{Tor(B,B_i)}$ is trivial. }

{\it Proof. } It is proven in [IM] (Proposition 2.6) that for any free resolution $P_{\bullet}$ of $B$ there are the following exact sequences(take $\Lambda=\Z$ and use a fact which states that ring $\Z$ has a global dimension one):

$$0 \to H_1( \varprojlim{P_{\bullet} \otimes B_i}) \to \varprojlim{Tor(B, B_i)} \to 0$$

So, we understood that $H_1( \varprojlim{P_{\bullet} \otimes B_i})$ is isomorphic to $\varprojlim{Tor(B, B_i)}$.

Let take $P_{\bullet}$ be a minimal resolution, i.e. such that $P_s=0$ when $s>1$. It is possible because global dimension of $\Z$ equals 1.

Let us note that the maps $P_s \otimes \varprojlim{B_i} \to  \varprojlim{(P_s \otimes B_i)}$ are embeddings because of Lemma 1 ($s=0,1$). 

Let consider a short exact sequence of complexes($C_{\bullet}$ is defined from this sequence): 
$$0 \to P_{\bullet}  \otimes (\varprojlim{B_i})  \to \varprojlim{(P_{\bullet} \otimes B_i)} \to C_{\bullet} \to 0 $$

Since $P_s$=0 for $s>1$, note  that $C_s=0$ for $s>1$. Thus we have the following exact sequence:

$$0 \to H_1(P_{\bullet}  \otimes (\varprojlim{B_i})) \to H_1(\varprojlim{(P_{\bullet} \otimes B_i)}) $$

But $H_1(P_{\bullet}  \otimes (\varprojlim{B_i}))=Tor(B, \varprojlim{B_i})$ by definition and we already have got that $H_1( \varprojlim{P_{\bullet} \otimes B_i})=\varprojlim{Tor(B, B_i)}$ so we are done. Q.E.D.

{\bf Corollary 5. }{\it $Tor (\varprojlim{A_i},\varprojlim{A_i}) \to \varprojlim {Tor(A_i,A_i)}$ is an embedding for any inverse system of abelian groups such that all maps $A_{i+1} \to A_i$ are epimorphic.}

{\it Proof. }  Let consider following commutative diagramm with $\psi$ being isomorphism:
 
\begin{tikzcd}
 Tor (\varprojlim{A_i},\varprojlim{A_j}) \arrow{r}{\varphi} \arrow[swap]{d}{f} & \varprojlim {Tor(A_i,A_i)} \arrow{d}{\psi} \\%
\varprojlim{Tor (A_i,\varprojlim{A_j})} \arrow{r}{g}& \varprojlim{\varprojlim{Tor (A_i,A_j)}}
\end{tikzcd}

But $f$ is monomorphism because of Statement 3 and $g$ is monomorphism because of Statement 3 and left exactness of $\varprojlim$. Q.E.D.

{\bf Statement 3. }{\it Let $B$ be a cotorsion group and $A$ be another abelian group. Let suppose that there are $i: A \to B$ and $\pi: B \to A$, such that $\pi i=n Id_A$ for $n \ge 1$. Then $A$ is a cotorsion group.  }

{\it Proof. } Since $Ext(\Q, B)=0$, a map $\pi i=nId_A$ induces zero endomorphism of $Ext(\Q, A)$. Let consider the following exact sequences(we denote $n$-torsion subgroup of $A$ by $A_n$):

$$0 \to A_n \to A \to nA \to 0$$

$$0 \to nA \to A \to A/nA \to 0$$

It is known that any $n$-torsion group is a cotorsion group (see [B], Theorem 8.5).

Thus the maps $A \to nA$ and $nA \to A$ induce isomorphisms on $Ext(\Q, -)$, and then so do $n Id_A$. But it induces zero endomorphism of $Ext(\Q, A)$, so we are done. Q.E.D.

{\it Proof of Theorem 3. } We already proved the first point of this theorem (Corollary 2), so let us prove second point of the theorem. It is known [Breen, section 6] that for any abelian group $A$ there is the following short exact sequence (here $L_1 \Lambda^2$ is first derived functor of functor $\Lambda^2$):

$$0 \to \Lambda^3(A) \to H_3(A) \to L_1 \Lambda^2(A) \to 0$$

Then for $A=\varprojlim{A_i}$ we have the following sequence:

$$0 \to \Lambda^3(\varprojlim{A_i}) \to H_3(\varprojlim{A_i}) \to L_1 \Lambda^2(\varprojlim{A_i}) \to 0$$

Since $\Lambda^3(\varprojlim{A_{i+1}}) \to \Lambda^3(\varprojlim{A_{i}})$ is an epimorphism for any $i$, ${\varprojlim}^1{\Lambda^3(\varprojlim{A_{i}})}=0$ and we have the following sequence:

$$0 \to \varprojlim{\Lambda^3(A_i)} \to \varprojlim{H_3(A_i)} \to \varprojlim{ L_1 \Lambda^2(A_i)} \to 0$$

Then $Ker[L_1 \Lambda^2(\varprojlim{A_i}) \to \varprojlim{ L_1 \Lambda^2(A_i)} ] \subseteq Ker [Tor (\varprojlim{A_i},\varprojlim{A_i}) \to \varprojlim {Tor(A_i,A_i)}]=0$ (here we used Corollary 5). Then, using Snake Lemma, we have the following short exact sequence:

$$0 \to Coker [\Lambda^3(\varprojlim{A_i} \to \varprojlim{\Lambda^3(A_i)}] \to Coker[H_3(\varprojlim{A_i}) \to \varprojlim{H_3(A_i)}] \to Coker[L_1 \Lambda^2(\varprojlim{A_i}) \to \varprojlim{ L_1 \Lambda^2(A_i)}] \to 0$$

It is obvious that any extension of a cotorsion group by a cotorsion group is a cotorsion group. And we already know that $Coker [\Lambda^3(\varprojlim{A_i} \to \varprojlim{\Lambda^3(A_i)}]$ is a cotorsion group. Then it is enough to prove that $Coker[L_1 \Lambda^2(\varprojlim{A_i}) \to \varprojlim{ L_1 \Lambda^2(A_i)}]$ is a cotorsion group.

It is known [Breen, sections 4 and 5] that for any abelian group $A$ there are maps $L_1\Lambda^2(A) \to Tor(A,A)$ and 
$Tor(A,A) \to L_1\Lambda^2(A)$, such that their composition equals to $2Id_{L_1\Lambda^2(A)}$. Then they induce natural maps $Coker[L_1 \Lambda^2(\varprojlim{A_i}) \to \varprojlim{ L_1 \Lambda^2(A_i)}] \to Coker[Tor (\varprojlim{A_i},\varprojlim{A_i}) \to \varprojlim {Tor(A_i,A_i)}]$ and $Coker[L_1 \Lambda^2(\varprojlim{A_i}) \to \varprojlim{ L_1 \Lambda^2(A_i)}] \to Coker[Tor (\varprojlim{A_i},\varprojlim{A_i}) \to \varprojlim {Tor(A_i,A_i)}]$, such that their composition is multiplication by two. Then, using Statement 3, we see that it is enough to prove that $Coker[Tor (\varprojlim{A_i},\varprojlim{A_i}) \to \varprojlim{Tor(A_i,A_i)}]$ is a cotorsion group. But if torsion subgroups of $A_i$ are groups of bounded exponent, then $Tor(A_i,A_i)$ are torsion groups of bounded exponent. Hence they are cotorsion groups [B, Theorem 8.5], and $\varprojlim {Tor(A_i,A_i)}$  is a cotorsion groups, because it is an inverse limit of cotorsion groups (it is obvious from our second definition of cotorsion groups that an inverse limit of cotorsion groups is itself a cotorsion group). Thus $Coker[Tor (\varprojlim{A_i},\varprojlim{A_i}) \to \varprojlim {Tor(A_i,A_i)}]$ is a cotorsion group as an image of a cotorsion group and we are done. Q.E.D.

We need the following statement in order to prove Theorem 2:

{\bf Statement 4. }{\it Let $A_i$ be an inverse system of torsion-free abelian groups indexed by $\mathbb{N}$ such that all maps $A_{i+1} \to A_i$ are epimorphisms and let $B$ be any torsion-free abelian group. Then the natural map $B \otimes \varprojlim{A_i} \to \varprojlim{(B \otimes A_i)}$ is an embedding. }

{\it Proof.} First let prove the statement for $B= \Q$. Let us consider the following short exact sequence:
$$0 \to \mathbb{Z} \to \Q \to \Q / \mathbb{Z} \to 0$$
After appliying the functor $- \otimes A_i$ for this sequence we get:
$$0 \to Tor(\Q / \mathbb{Z}, A_i)\to A_i \to \Q \otimes A_i \to \Q / \mathbb{Z} \otimes A_i \to 0$$
Since ${Tor(\Q / \mathbb{Z}, A)}=t(A)$ for any abelian group $A$, we have:
$$0 \to  t(A_i) \to A_i \to \Q \otimes A_i \to \Q / \mathbb{Z} \otimes A_i \to 0$$
Then, using left-exactness of $\varprojlim$, we have:
$$0 \to  \varprojlim{t(A_i)} \to \varprojlim{A_i} \to \varprojlim{(\Q \otimes A_i)}$$
Let us apply the exact functor $\Q \otimes -$ for this sequence:
$$0 \to  \Q \otimes \varprojlim{t(A_i)} \to \Q \otimes \varprojlim{A_i} \to \Q \otimes \varprojlim{(\Q \otimes A_i)}$$
Let us note that $\Q \otimes \varprojlim{(\Q \otimes A_i)} \cong \varprojlim{(\Q \otimes A_i)}$, because $\varprojlim{(\Q \otimes A_i)}$ is a $\Q$-vector space. It means that we got the following:
$$0 \to  \Q \otimes \varprojlim{t(A_i)} \to \Q \otimes \varprojlim{A_i} \to \varprojlim{(\Q \otimes A_i)}$$
But in the our case $t(A_i)=0$ for any $i$. So we have:
$$0 \to \Q \otimes \varprojlim{A_i} \to \varprojlim{(\Q \otimes A_i)}$$

Now let us prove the statement in the case $B= {\Q}^{\oplus I}$, where $I$ is any cardinal. Using the sequence for $B=\Q$  we get:
$$0 \to  {(\Q \otimes \varprojlim{t(A_i)})}^{\oplus I} \to {(\Q \otimes \varprojlim{A_i})}^{\oplus I} \to {(\varprojlim{(\Q \otimes A_i)}) }^{\oplus I} $$
Since it is obvious that the map ${(\varprojlim{(\Q \otimes A_i)}) }^{\oplus I} \to {\varprojlim{((\Q \otimes A_i)^{\oplus I})} }$ is injective, we have:
$$0 \to  {(\Q \otimes \varprojlim{t(A_i)})}^{\oplus I} \to {(\Q \otimes \varprojlim{A_i})}^{\oplus I} \to  {\varprojlim{((\Q \otimes A_i)^{\oplus I})} }$$
But in the our case $t(A_i)=0$ for any $i$. So we have the exactness of the following sequence:
$$0 \to {(\Q \otimes \varprojlim{A_i})}^{\oplus I} \to  {\varprojlim{((\Q \otimes A_i)^{\oplus I})} }$$
Since the functor $- \otimes \A$ commutes with direct sums for any $A$, we proved the statement in the case when $B$ is any $\Q$-vector space. 

Now let us prove the statement for any torsion-free $B$. 
Let $Q_B$ be the injective hull of $B$. Since $B$ is torsion-free, $Q_B=\Q^{\oplus I}$ for some cardinal $I$.
Let us consider the following short exact sequence($B'$ is defined from this sequence):
$$0 \to B \to Q_B \to B' \to 0$$
Since $\varprojlim{A_i}$ is torsion-free and ${\varprojlim}^1{B \otimes A_i}=0$ because the maps $B \otimes A_{i+1} \to B \otimes A_i$ are epimorphisms, it gives us two sequences:
$$0 \to B \otimes \varprojlim{A_i} \to Q_B \otimes \varprojlim{A_i} \to B' \otimes \varprojlim{A_i}  \to 0$$
$$0 \to \varprojlim{(B \otimes A_i)} \to \varprojlim{(Q_B \otimes A_i)} \to \varprojlim{(B' \otimes A_i)}  \to 0$$
Then, using Snake Lemma for this two sequences and natural maps between their elements, we have:
$$0 \to Ker[B \otimes \varprojlim{A_i} \to \varprojlim{(B \otimes A_i)}] \to Ker [Q_B \otimes \varprojlim{A_i} \to \varprojlim{(Q_B \otimes A_i)}]$$
Since $Q_B$ is a $\Q$-vector space, we already proved that $Ker [Q_B \otimes \varprojlim{A_i} \to \varprojlim{(Q_B \otimes A_i)}]=0$. Then $Ker[B \otimes \varprojlim{A_i} \to \varprojlim{(B \otimes A_i)}]=0$ and we are done. Q.E.D. 

{\it Proof of Theorem 2. }Let consider following commutative diagram with $\psi$ being isomorphism.

\begin{tikzcd}
 \varprojlim{A_i} \otimes \varprojlim{A_j} \arrow{r}{\varphi} \arrow[swap]{d}{f} & \varprojlim {(A_i \otimes A_i)} \arrow{d}{\psi} \\%
\varprojlim{(A_i \otimes \varprojlim{A_j})} \arrow{r}{g}& \varprojlim{\varprojlim{ (A_i \otimes A_j)}}
\end{tikzcd}

Then we see that $f$ is embedding because of Statement 4 and $g$ is embeddging because of Statement 4 and left-exactness of $\varprojlim$. This implies that $\phi$ is also embedding. Then, since there is a natural embedding $\Lambda^2(A) \to A \otimes A$ and $H_2$ is naturally isomorphic to $\Lambda^2$ [Breen, section 6]  for abelian groups, we proved the theorem for $n=2$.

Now let us note that $H_n$ is naturally isomorphic to $\Lambda^n$ on the category of torsion-free abelian groups also for any $n>2$ [Breen, p. 214, (1.12)]. 

Let us prove the theorem by induction on $n$.

We can assume that $n \ge 3$. Let us note that the map $({\varprojlim{A_i}})^{\otimes n}   \to \varprojlim{{(A_i)}^{\otimes n}}$ may be decomposed up to isomorphism in the following way:

$({\varprojlim{A_i}})^{\otimes n} \cong \varprojlim{A_i} \otimes (\varprojlim{A_j})^{\otimes (n-1)}  \to \varprojlim{A_i} \otimes  \varprojlim{{(A_j)}^{\otimes (n-1)}} \to \varprojlim_i {\varprojlim_j{ A_i \otimes  {(A_j)}^{\otimes (n-1)}}} \cong \varprojlim{{(A_i)}^{\otimes n}}$.  But all maps in the decomposition are monic because of inductional assumption, left-exactness of $\varprojlim$ and exactness of $- \otimes B$ for torsion-free $B$. Then the map $({\varprojlim{A_i}})^{\otimes n}   \to \varprojlim{{(A_i)}^{\otimes n}}$ is also monic, and then so is $H_n(\varprojlim{A_i}) \to \varprojlim{H_n(A_i)}$. Q.E.D.

\section{Applications to topology}

{\bf Theorem 4. }{\it Let $X_i$ be an inverse system  of pointed connected spaces, such that all maps $X_{i+1} \to X_i$ are Serre fibrations, all $\pi_1(X_i)$ are abelian, all maps $\pi_1(X_{i+1}) \to \pi_1(X_{i})$ are epimorphisms and  $\pi_2(X_i)_{i \in \mathbb{N}}$ satisfy the Mittag-Leffler condition. Then:

(1) Cokernel of the natural map $H_2(\varprojlim{X_i}) \to \varprojlim {H_2(X_i)}$ is a cotorsion group.

(2) Suppose, additionally, that all $\pi_1(X_i)$ are torsion-free. Then kernel of the natural map $H_2(\varprojlim{X_i}) \to \varprojlim {H_2(X_i)}$ is a cotorsion group.

(3) Suppose that condition (2) is satisfied  and, additionally, $\pi_3(X_i)_{i \in \mathbb{N}}$ satisfy the Mittag-Leffler condition. Then the natural map $H_2(\varprojlim{X_i}) \to \varprojlim {H_2(X_i)}$ is embedding.}

{\it Proof. }Let us denote $\Pi_2(X):=Im[\pi_2(X) \to H_2(X)]$. It is obvious that $\Pi_2$ is a functor. Then we have the following sequence:
$$0 \to \Pi_2(X) \to H_2(X) \to H_2(\pi_1(X)) \to 0$$
Then we get the following sequences, feeding $X=\varprojlim{X_i}$ and $X=X_i$:
$$0 \to \Pi_2(X_i) \to H_2(X_i) \to H_2(\pi_1(X_i)) \to 0$$
It gives us the following sequence:
$$0 \to \varprojlim{\Pi_2(X_i)} \to \varprojlim{H_2(X_i)} \to \varprojlim{H_2(\pi_1(X_i))} \to {\varprojlim}^1{\Pi_2(X_i)} \to {\varprojlim}^1{H_2(X_i)} \to {\varprojlim}^1{H_2(\pi_1(X_i))} \to 0$$

Since $\pi_2(X_i) \to \Pi_2(X_i)$ is epimorphism, ${\varprojlim}^1{\pi_2(X_i)} \to {\varprojlim}^1{\Pi_2(X_i)}$ it is well-known (e.g. see [Brown, p. 42, Theorem 5.2]) that for any space $X$ there is the following exact sequence: 
$$\pi_2(X) \to H_2(X) \to H_2(\pi_1(X)) \to 0$$. also epimorphism. Then ${\varprojlim}^1{\pi_2(X_i)}= {\varprojlim}^1{\Pi_2(X_i)}=0$, because $\pi_2(X_i)$ satisfy the Mittag-Leffler condition. So we have the following sequence:
$$0 \to \varprojlim{\Pi_2(X_i)} \to \varprojlim{H_2(X_i)} \to \varprojlim{H_2(\pi_1(X_i))} \to 0$$
Also we have the following sequence:
$$0 \to \Pi_2(\varprojlim{X_i}) \to H_2(\varprojlim{X_i}) \to H_2(\pi_1(\varprojlim{X_i})) \to 0$$
We have the following sequences, because ${X_i}_{i \in \mathbb{N}}$ is a tower of pointed fibrations [GJ, VI Proposition 2.15]:
$$0 \to {\varprojlim}^1{\pi_2(X_i)} \to \pi_1(\varprojlim{X_i}) \to \varprojlim{\pi_1(X_i)} \to 0$$
Since $\pi_2(X_i)$ satisfy the Mittag-Leffler condition, we have:
$$\pi_1(\varprojlim{X_i}) \cong \varprojlim{\pi_1(X_i)}$$
Finally, we have two sequences:
$$0 \to \varprojlim{\Pi_2(X_i)} \to \varprojlim{H_2(X_i)} \to \varprojlim{H_2(\pi_1(X_i))} \to 0$$
$$0 \to \Pi_2(\varprojlim{X_i}) \to H_2(\varprojlim{X_i}) \to H_2(\varprojlim{\pi_1(X_i)}) \to 0$$.

Let us note that $Coker[\varprojlim{\Pi_2(X_i)} \to \Pi_2(\varprojlim{X_i})]=0$ and $Ker[\varprojlim{\Pi_2(X_i)} \to \Pi_2(\varprojlim{X_i})] \cong {\varprojlim}^1{\pi_3(X_i)} $ . It follows from the following sequence:
$$0 \to {\varprojlim}^1{\pi_3(X_i)} \to \pi_2(\varprojlim{X_i}) \to \varprojlim{\pi_2(X_i)} \to 0$$
Then, using Snake Lemma for this two sequences and the natural maps between them, we have:
$$0 \to {\varprojlim}^1{\pi_3(X_i)} \to Ker[\varprojlim{H_2(X_i)} \to H_2(\varprojlim{X_i})] \to Ker[\varprojlim{H_2(\pi_1(X_i))} \to H_2(\varprojlim{\pi_1(X_i)})] \to 0$$ and $$Coker[H_2(\varprojlim{X_i}) \to \varprojlim{H_2(X_i)}] \cong Coker[H_2(\varprojlim{\pi_1(X_i)}) \to \varprojlim{H_2(\varprojlim{(X_i))}}$$
Then point (1) of the Theorem follows from point (1) of Theorem 1. Let us prove points (2) and (3) of the Theorem. It follows from Theorem 2 for $n=2$ that $Ker[\varprojlim{H_2(\pi_1(X_i))} \to H_2(\varprojlim{\pi_1(X_i)})]=0$. Then we have:
$$Ker[\varprojlim{H_2(X_i)} \to H_2(\varprojlim{X_i})] \cong {\varprojlim}^1{\pi_3(X_i)} $$ 
So we proved point (3) and point (2). Q.E.D.

{\bf Theorem 5. }{\it Let $Y_i$ be a sequence  of spaces. Then cokernel of the map $H_k(\displaystyle \prod_{i=1}^{\infty} Y_i) \to \varprojlim_{n}{H_k(\prod_{i=1}^{n} Y_i)}$ is a cotorsion group for every $k$.}

{\it Proof. } Since composition of the natural maps $H_k(\displaystyle \prod_{i=1}^{n} Y_i) \to H_k(\displaystyle \prod_{i=1}^{\infty} Y_i)$ and $H_k(\displaystyle \prod_{i=1}^{\infty} Y_i) \to H_k(\displaystyle \prod_{i=1}^{n} Y_i)$ is the identity map, the map $H_k(\displaystyle \prod_{i=1}^{\infty} Y_i) \to H_k(\displaystyle \prod_{i=1}^{n} Y_i)$ is surjective. So we are done because of Statement 1. Q.E.D.

{\it Remark. }It is easy to see from this from this proof that the same fact holds for any category which has infinite products instead of category of spaces and for and functor $F$ from this category to $Ab$ instead of $H_k$. But I formulated it in this way because it seems more natural in this section.

We will also prove the following Theorem, which does not follow from our previous results but follows from Shelah's and Barnea's. It is connected with Theorem 3, so I found quite natural to formulate and prove it here. This Theorem is a generalization of Shelah's and Barnea's [Corollary 0.0.9. BS]. Actually we show that it is not necessary to assume that $\pi_2(X_{i})$ satisfies the Mittag-Leffler conditon. 

{\bf Theorem 6. }{\it Let $X_i$ be an inverse system  of pointed connected spaces, such that all maps $X_{i+1} \to X_i$ are Serre fibrations and $\pi_1(X_{i})$ satisfies the Mittag-Leffler condition. Then the natural map $H_1(\varprojlim{X_i}) \to \varprojlim {H_1(X_i)}$ is surjective and its kernel is a cotorsion group.}

{\it Proof.} Let us consider the following exact sequence:
$$0 \to {\varprojlim}^1{\pi_1(X_i)} \to \pi_0(\varprojlim{X_i}) \to \varprojlim{\pi_0(X_i)} \to 0$$
Then $\varprojlim{X_i}$ is also connected and we have:
$$H_1(\varprojlim{X_i}) \cong {\pi_1(\varprojlim{X_i})}_{ab}$$
It follows that: $Ker[H_1(\varprojlim{X_i}) \to \varprojlim {H_1(X_i)}]=Ker[{(\pi_1(\varprojlim{X_i}))}_{ab} \to \varprojlim{{(\pi_1(X_i))}_{ab}}]$ and $Coker[H_1(\varprojlim{X_i}) \to \varprojlim {H_1(X_i)}]=Coker[{(\pi_1(\varprojlim{X_i}))}_{ab} \to \varprojlim{{(\pi_1(X_i))}_{ab}}]$.
Also we have the following sequence:
$$0 \to {\varprojlim}^1{\pi_2(X_i)} \to \pi_1(\varprojlim{X_i}) \to \varprojlim{\pi_1(X_i)} \to 0$$
Thus we have the following exact sequence from 5-term exact sequence:
$${\varprojlim}^1{\pi_2(X_i)} \to {(\pi_1(\varprojlim{X_i}))}_{ab} \to {(\varprojlim{\pi_1(X_i)})}_{ab} \to 0$$
Then, since ${\varprojlim}^1{\pi_2(X_i)}$ is a cotorsion group by the Theorem of Huber(Theorem 1, [H]) and  any quotient group of a cotorsion group is itself a cotorsion group, $Ker[{(\pi_1(\varprojlim{X_i}))}_{ab} \to {(\varprojlim{\pi_1(X_i)})}_{ab}]$ is a cotorsion group. Now consider the map ${(\varprojlim{\pi_1(X_i)})}_{ab} \to \varprojlim{{(\pi_1(X_i))}_{ab}}$. It is surjective and its kernel is a cotorsion group by Theorem 0.0.1 of Shelah and Barnea from [BS]. Now note that we can decompose the map ${(\pi_1(\varprojlim{X_i}))}_{ab} \to \varprojlim{{(\pi_1(X_i))}_{ab}}$  in the following way:
$${(\pi_1(\varprojlim{X_i}))}_{ab} \to (\varprojlim {\pi_1({X_i}))}_{ab} \to  \varprojlim{{(\pi_1(X_i))}_{ab}}$$
Then the map is surjective as a composition of two surjective maps and its kernel is an extension of a cotorsion group by a cotorsion group. Thus it is obvious from our second definition of cotorsion groups that the kernel itself is a cotorsion group . Q.E.D.

\section{References}

[B] Reinhold Baer. The subgroups of elements of finite order of an Abelian group, Ann. of Math., 37, (1936), 766-781.

[Bous] A. K. Bousfield. Homological localization towers for groups and $\pi$-modules, Mem. Amer. Math. Soc, no. 186, 1977.

[Brown] K.S. Brown, Cohomology of Groups, Graduate Texts in Mathematics, vol. 87, Springer, Berlin, 1982.

[Breen] Lawrence Breen. On the functorial homology of abelian groups, Journal of Pure and Applied Algebra 142 (1999) 199-237

[BS] Ilan Barnea, Saharon Shelah. The abelianization of inverse limits of groups, Israel Journal of Mathematics, August 2018, Volume 227, Issue 1, pp. 455-483.

[D] W.G. Dwyer. Homological localization of $\pi$-modules, Journal of Pure and Applied Algebra, Volume 10, Issue 2, November 1977, Pages 135-151.

[F] L. Fuchs. Infinite Abelian Groups, Academic Press, Jan 1, 1970.

[GJ] Goerss P. G., Jardine J. F. Simplicial Homotopy Theory, Progress in Mathematics, Vol. 174, Birkhauser,
Basel, 1999.

[H] Warfield, R. B., Jr., Huber, Martin. On the values of the functor ${\varprojlim}^1$, Arch. Math. (Basel) 33 (1979/80), no. 5, 430-436. 

[Hill] Jonathan Hillman. Algebraic Invariants of Links, Series on Knots and Everything, Vol. 52.

[IM] Sergei O. Ivanov, Roman Mikhailov. On a problem of Bousfield for metabelian groups, Advances in Mathematics, Volume 290, 26 February 2016, Pages 552-589.

[IM2] Sergei O. Ivanov, Roman Mikhailov. A finite Q-bad space, arXiv:1708.00282, 1 August 2017 

[IM3] Sergei O. Ivanov, Roman Mikhailov. On discrete homology of a free pro-$p$-group, Compositio Mathematica, Volume 154, Issue 10, October 2018 , pp. 2195-2204.

[KM] O. Kharlampovich, A. Miasnikov. Implicit function theorem over free groups and genus problem Proceedings of the Birmanfest, AMS/IP Studies in Advanced Mathematics, v. 24, 2001, 77-83.

[MP] Roman Mikhailov, Inderbir Singh Passi. Lower Central and Dimension Series of Groups, Lecture Notes in Mathematics, 2009.

[W] Charles A. Weibel. An Introduction to Homological Algebra, Cambridge University Press, 1994.

\bigskip

Laboratory of Modern Algebra and Applications, St. Petersburg State University, 14th Line, 29b, Saint Petersburg, 199178, Russia 

{\it E-mail address}: akhtyamoff1997@gmail.com

\end{document}